\documentclass[a4paper,oneside,11pt]{article}
\date{}

\usepackage[latin1]{inputenc}
\usepackage[T1]{fontenc}
\usepackage [paper=a4paper, 
	inner=2.5cm, 
	outer=2.5cm, 
	top=2.8cm, 
	bottom=2.8cm,]{geometry}

\usepackage{amsmath,amssymb,amsthm,amsfonts}
\usepackage{indentfirst}
\usepackage{graphicx}
\usepackage{enumitem}
\usepackage{verbatim}
\usepackage{mathrsfs}
\usepackage{tikz}
\usepackage[colorlinks,allcolors=red]{hyperref}
\usepackage{url}
\usepackage{bbm}
\usepackage[toc,page]{appendix}
\usepackage{epigraph}
\usetikzlibrary{matrix}
\usetikzlibrary{arrows}
\usetikzlibrary{graphs}
\usetikzlibrary[graphs]
\def\la{\big\langle}
\def\ra{\big\rangle}

\def\ds{\displaystyle}
\def\forall{\hbox{for all}~}
\def\L{{\bf L}}
\def\D{{\cal D}}

\def\argmin{\hbox{arg}\!\min}
\def\bfv{{\bf v}}

\def\bfw{{\bf w}}
\def\bfu{{\bf u}}

\def\bfn{{\bf n}}
\def\bfe{{\bf e}}

\def\bfk{{\bf k}}
\def\ve{\varepsilon}

\def\A{{\cal A}}

\def\R{{\mathbb R}}
\def\N{{\mathbb N}}

\def\implies{\Longrightarrow}

\def\v{\vskip 1em}

\def\begi{\begin{itemize}}
\def\endi{\end{itemize}}
\def\C{{\cal C}}
\def\wto{\rightharpoonup}

\def\ov{\overline}

\def\bega{\begin{array}}
\def\enda{\end{array}}

\def\bel{\begin{equation}\label}
\def\eeq{\end{equation}}
\def\sqr#1#2{\vbox{\hrule height .#2pt
\hbox{\vrule width .#2pt height #1pt \kern #1pt
\vrule width .#2pt}\hrule height .#2pt }}
\def\square{\sqr74}
\def\endproof{\hphantom{MM}\hfill\llap{$\square$}\goodbreak}

\newtheorem{theorem}{Theorem}[section]
\newtheorem{lemma}[theorem]{Lemma}

\newtheorem{definition}[theorem]{Definition}
\newtheorem{remark}[theorem]{Remark}

\theoremstyle{plain}

\newtheorem*{thm*}{Theorem}
\newtheorem{teor*}{Theorem}

\theoremstyle{definition}

\theoremstyle{remark}

\numberwithin{equation}{section}

\title{Soil searching by an artificial root} 

\author{ Fabio Ancona$^{(1)}$, Alberto Bressan$^{(2)}$, Maria Teresa Chiri$^{(2)}$\\
\\
 {\small $^{(1)}$ Dipartimento di Matematica ``Tullio Levi-Civita'', Universit\`a di Padova, }\\  
   {\small $^{(2)}$ Department of Mathematics, Penn State University}\\ 
 \\  {\small E-mails: ~ancona@math.unipd.it,  ~axb62@psu.edu, ~ mxc6028@psu.edu}
}
\begin{document}

\maketitle
\begin{abstract}
We model an artificial root which grows in the soil for underground prospecting. Its evolution is described by a  controlled 
system of two integro-partial differential equations: 
one for the growth of the body and the other for the 
 elongation 
of the tip. 
At any given time, the angular velocity of the root is obtained by solving a minimization problem with state constraints. We prove the 
existence of solutions to the evolution problem, up to the first time where 
a  "breakdown configuration" is reached. Some numerical simulations are
performed, to test the effectiveness of our  feedback control algorithm.
\end{abstract}

\section{Introduction}
We consider a mathematical model for   a robotic root,
which is able to move inside a medium by growing and bending.
These artificial  roots are made of a tubular body, a growing head
(realized by depositing new material adjacent to the tip)
and a sensorized tip that commands the robotic behaviors. 
These mechanisms are inspired by the so-called {\it apical growth} or {\it elongation from the tip}
 process  in plant roots: 
newly generated cells by mitosis  
move from the meristematic region near the root apex to the elongation region
where they expand axially, while mature cells of the root remain stationary and in contact with the soil.
Root-like devices, made of deformable materials
(continuum  soft robots), are particularly efficient to act in a 
space-constrained
and unstructured environment 
(e.g. exploring uncertain terrains searching for life or resources in the soil, 
accomplishing predetermined tasks in environments with solid obstacles~\cite{8460777, 6496902},  or performing minimally invasive surgery~\cite{burgner})
and to interact more safely with humans (e.g. rescuing people trapped by 
earthquakes or avalanches,
in hazardous areas).

Within the growing interest in biologically inspired technologies~\cite{pfeifer}, 
robots inspired by the stems and roots of plants are receiving a distinctive attention because plants 
are the most efficient and 
compliant living organisms for soil exploration.
Therefore, their remarkable features can be exploited in artificial systems to
ensure energy efficiency and minimize friction while exploring and searching
( see (\cite{Okamura2, Mazzolai1, Mazzolai2, robotics7030058}).

In the present paper we provide a mathematical analysis of two simple models of continuum robotic roots with apical growth,
which explore the soil in order to reach a certain target. 
This study should help in developing better control strategies and optimizing the implementation of these new technologies. 

In the same spirit of the models introduced in~\cite{07a9e864592f44a99f73752cdcb6e0a8,
MR3650339} for the growth of a plant stem or a vine,
at each time $t$ the robotic root will be described 
by a curve $\gamma(t,\cdot)$ in 3-dimensional space.
The model takes into account the fact that the artificial root is able to bend when the soil is too hard or when it finds an obstacle, and it can even  contract when it is not able to outflank an hindrance. 
The evolution of the robotic root is governed by a system of two partial differential equations, one for the growth of the body, and the other one for the extension of the tip. The angular velocity of the root is obtained as solution of an optimization problem with state constraint.
In fact, it is
physically meaningful that the robot bends minimizing the 
instantaneous deformation energy and the cost to move sand or to drill in high-density regions.
\v


For a general description of plant and root development from a biological point of view we refer to~\cite{leyser}
and~\cite{Dexter1987}, respectively,
while for a mathematical theory of biological growth we refer to~\cite{MR2014211, MR2288420,1686075}.
\v

The remainder of the paper is organized as follows: in Section~\ref{S31} we introduce the basic models for the evolution of our artificial root. 
In Sections~\ref{S32} and \ref{S33} we describe respectively a rigid and an elastic model for the robot, introducing a restarting procedure when it is no longer 
able to grow because of an obstacle. 
Most of our analysis is focused on the elastic model, since it is the most realistic, and it provides more interesting theoretical challenges.  
This is covered in Section~\ref{S34}. 
Section~\ref{S35} is devoted to the construction of a  solution to the evolution problem.  Finally, in Section~\ref{S36} we present some numerical
simulations. 

 \section{Modeling the growth of an artificial root}\label{S31}
\label{s:0}
\setcounter{equation}{0}

We seek to model the movement of an artificial root, 
which penetrates the ground searching for water, chemical compounds, 
cavities where earthquake survivors may be trapped, etc$\ldots$

For this purpose, we need to assign:
\begi
\item[{\bf (I)}] An equation describing the scalar 
velocity at which the length of the root grows.
\item[{\bf (II)}]  An equation describing how the orientation $\bfk(t)$ 
of the tip of the root changes in time.
\item[{\bf (III)}] A rule determining when the growth stops, the root shrinks to an earlier
configuration, and then starts growing in a new direction.
\endi

 \subsection{A rigid root}\label{S32}
\label{s:1}
\setcounter{equation}{0}
In this first rigid model we assume that the root grows and bends from the tip but the portion 
 of the root that 
has already grown does not move. 
  In this case, apart from restarting times,
the growth can be determined by a second order ODE.

\v
{\bf 1.} At any time $t$, the root will be described by a $\C^{1,1}$ curve
with values in $\R^3$:
\bel{ga}s~\mapsto~ \gamma(t,s),\qquad s\in [0, t],\eeq
parameterized by arc-length.  
The point $ \gamma(t,s)$ denotes the position 
at time $s\in [0,t]$ of the tip of the root grown until time $t$.
Since we are assuming that  when the root apex grows no deformation of the rest of the body
is undertaken,
we have
\bel{rigid-growth}
\gamma_t(t,s)=0\qquad  \forall \ s\in [0,t]\,.
\eeq
For simplicity we assume that the root grows at unit speed, i.e. that $|\gamma_s(t,s)|=1$
for all $s\in [0.t]$, so that at time $t$
its length is simply $\ell(t)=t$.
We denote by
\bel{tv}\bfk(t,s)~\doteq~\gamma_s(t,s)\eeq
the unit tangent vector to the curve.
\v
{\bf 2.} 
For  convenience we shall denote by
\bel{Pk}P(t)~\doteq~\gamma(t,t),\qquad\qquad \bfk(t)~\doteq~\bfk(t,t)~=~\gamma_s(t,t)\eeq
respectively  the position and the orientation of the root tip.

Two types of control will be implemented. Namely:
\begi
\item[(i)]  At every time $t$, the orientation of the tip is modified.
\item[(ii)] At a finite number of times $0<t_1<t_2<\cdots$, if the root hits an obstacle, a restarting
procedure is adopted and growth is restarted in a different direction.
\endi

\v
{\bf 3.} 
At any time $t$, we assume that the  direction $\bfk(t)$ of the root tip can change, subject to a 
maximum possible curvature.
Namely, 
\bel{g4} \dot \bfk(t)~=~\bfu(t)\times \bfk(t),\qquad\qquad |\bfu(t)|~\leq~\kappa_0\eeq
Here $\bfu(\cdot)$ is a measurable control function, describing the angular velocity at which 
we bend the tip of the root.
\v
{\bf 4.} We now introduce a feedback rule, assigning the control $\bfu$ in terms of
the position of the tip. 
Toward this goal, we consider two scalar functions.
\begi
\item A function $\phi:\R^3\mapsto \R$, given at the beginning and never modified afterwards. 
Roughly speaking one may think of $\phi(x)$   as the expected distance from a (random) target.
\item A function $\psi:\R^3\mapsto \R$, measuring a kind of closeness to previously reached points.
A bit more precisely, if $\Gamma(t)$ is the union of all points reached by the artificial root 
during previously failed attempts, we could define
\bel{G}
\psi(x)~\doteq~\int_{y\in \Gamma(t)} e^{-|x-y|} \, dy,\eeq
\endi

In view of~\eqref{rigid-growth}, \eqref{Pk}, \eqref{g4},  the growth of the root can then be described by 
\bel{rg} 
\dot P(t)~=~\bfk(t),
\qquad\qquad \dot \bfk(t) ~=~\bfu(t)\times \bfk(t),\eeq
where the control vector $\bfu$ is determined by
\bel{uopt}\bfu(t, P(t))~=~\argmin_{|\omega|\leq \kappa_0} \Big\langle \omega\times \bfk(t)\,,\,
\nabla \phi (P(t)) + \nabla \psi(P(t))\Big\rangle.\eeq
Note: in (\ref{G}), instead of $e^{-|x-y|}$, one could use some other kernel
$K(|x-y|)$, where $s\mapsto K(s)$ is a smooth, decreasing function.

 In some cases, the control in (\ref{uopt}) may not be unique.
To achieve a well posed evolution equation, it is convenient to replace (\ref{uopt})
with 
\bel{uopt1}\bfu(t,P(t))~=~\argmin_{|\omega|\leq \kappa_0}\Big\{
 \la \omega\times \bfk(t)\,,\nabla \phi (P(t))+ \nabla \psi(P(t))\ra + \ve |\omega|^2\Big\}.\eeq
The strict convexity of the right hand side yields a unique minimizer, depending 
Lipschitz continuously on $P$.  Of course, other approximations are possible, with 
$\C^\infty$ dependence on  $\bfk, \nabla\phi, \nabla\psi$.

Having done this, the local existence and uniqueness of solutions to the evolution problem~\eqref{rg} is trivial.   




\begin{figure}[ht]
\centerline{\hbox{\includegraphics[width=13cm]{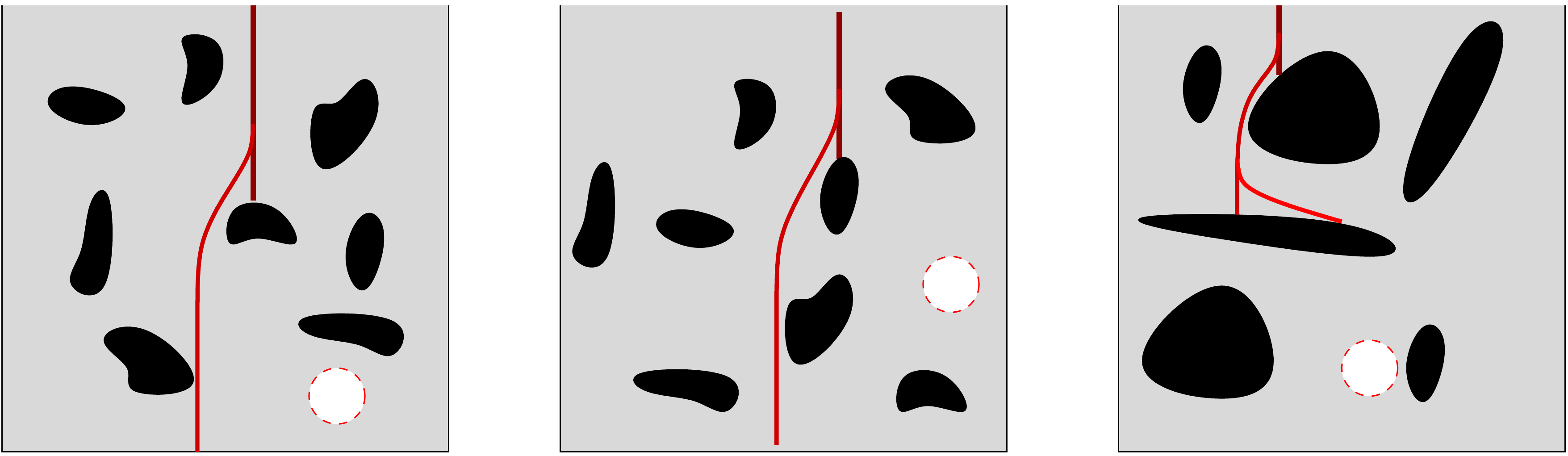}}}
\caption{\small   Rocks with different average sizes should 
determine different restarting strategies.  The white circles denote the targets.}
\label{f:ar4}
\end{figure}

{\bf 5.} In addition, there must be a rule prescribing when the growth should
stop and the root should go back to a previous configuration, and restart.

We denote by $\gamma_i(\cdot,\cdot)$ the curve constructed at the $i$-th trial run,
and write $P_i(t) = \gamma_i(t,t)$  for the position of the tip of this curve at time $t$.
The $i$-th curve starts with length $t_i^-$ and grows up to 
a maximum length $t_i^+$.   

Here $h:\R^3\mapsto [0,1]$ is some scalar function, not known a priori,
accounting for the ``hardness" 
of the soil.   For example, $h=0$ 
when the soil offers no resistance (say, inside a crack), while $h>\!>1$ along an impenetrable obstacle
(say, a wall).

A reasonable stopping rule is then
\bel{tide}
t_i^+~\doteq~\max\Big\{ t > t_i^-\,;~~h(P(t)) \leq h_0\Big\}\,.\eeq
Here $h_0>0$ is a threshold value, assigned a priori.
The idea is that, if the tip of the root encounters very hard soil or a rock,
it is better go back and restart in another direction.

The restarting procedure must also be carefully defined.  
Calling $t_i^+$ the length of the root when we decide to restart,
we choose a new length
$t_{i+1}^-<t_i^+$ and set 
\bel{gi+}
\gamma(s,t_{i+1}^-)~=~\gamma(s,t_i^+)\qquad \qquad s\in [0, t_{i+1}^-].\eeq
In general, the new length will be determined by a restarting function:
\bel{R} t_{i+1}^-~=~R\Big(t_i^+, t_i^-, \psi(P(t_i^+))\Big).\eeq
By (\ref{R}), we are saying that the amount by which the root shrinks
depends on 
\begi
\item the length $t_i^+$ and the change in the length $t_i^+ -t_i^-$, achieved during the previous failed attempt,
\item the density $\psi(P(t_i^+))$, measuring how much the region near  $P(t_i^+)$ has already been explored. 
\endi

A first basic example of restarting function that one can consider is 
\bel{R1} R_1\Big(t_i^+, t_i^-, \psi(P(t_i^+))\Big) ~=~ t_i^+ + c\Big( e^{ -\psi(P(t_i^+))}-1\Big) \Big(t_i^+ - t_i^-\Big)\eeq
with $c>1$. This can be interpreted as follows: if the value 
of $\psi (P(t_i^+))$ is large (which means that $P(t_i^+)$ is close to many trajectories already explored in previous attempts), then we need  to go back a strictly greater quantity than the last change of length $t_i^+ -t_i^- $, if  $\psi (P(t_i^+))$ is small, we do not need to restart so far from the last starting point $P(t_i^-)$.\\
As an alternative restarting function, we propose:
\bel{R2} R_2\Big(t_i^+, t_i^-, \psi(P(t_i^+))\Big) ~=~ t_i^+ + c\Big( 1+ 2\rho \chi_{[0,\rho]}(|P(t_i^+)-P(t_i^-)|)  \Big) \Big( e^{ -\psi(P(t_i^+))}-1\Big) \Big(t_i^+ - t_i^-\Big).\eeq
Compared with (\ref{R1}), we have now inserted a factor measuring the length of the last attempt. This is motivated by the fact that, if the last elongation is small (in the interval $[0,\rho]$, for some fixed threshold value $\rho>0$), then probably the root has not enough space to move, therefore it is reasonable to retreat back 
by a length much greater then $t_i^+ -t_i^- $.\\
\begin{figure}[ht]
\centerline{\hbox{\includegraphics[width=6cm]{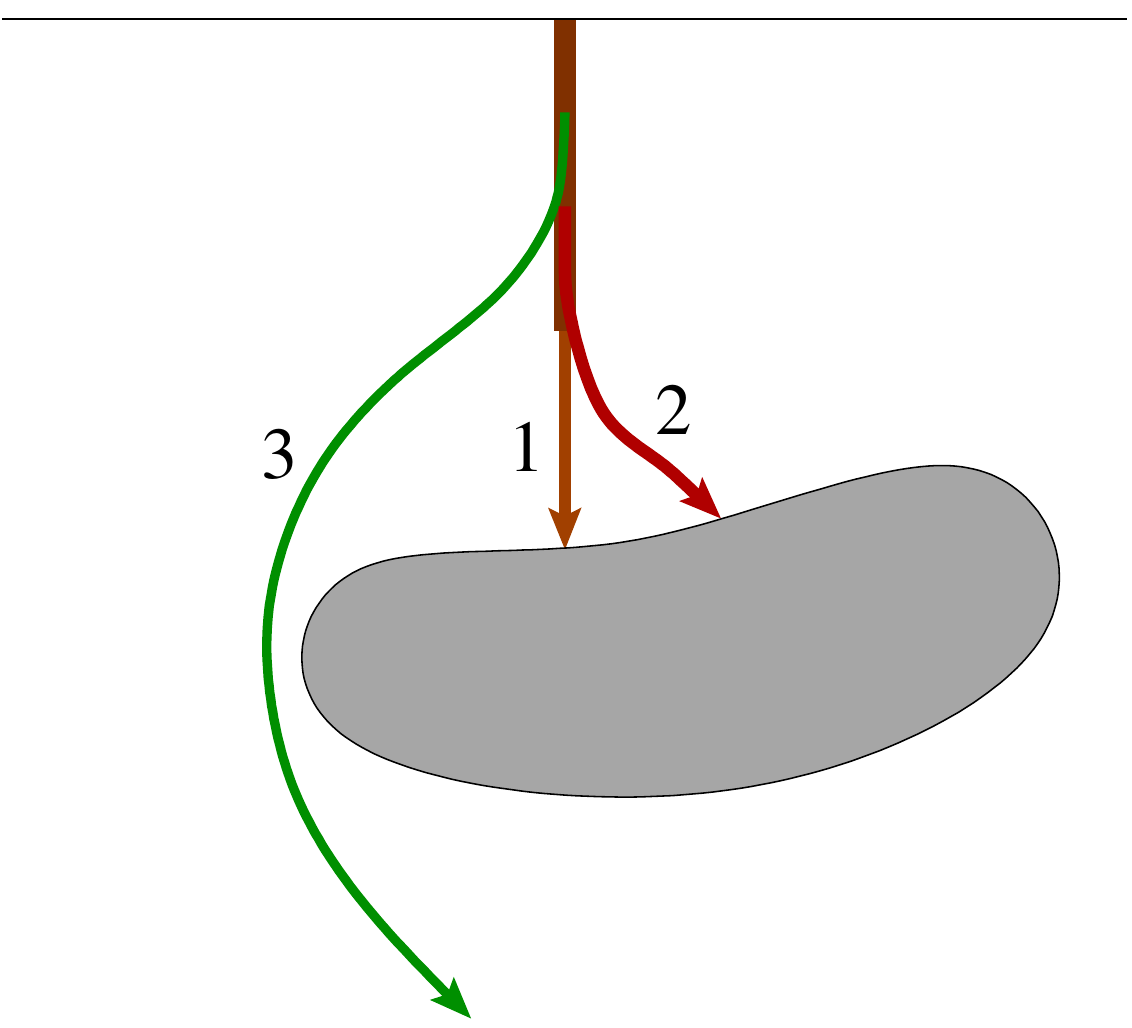}}}
\caption{\small   Three attempts at growing a root past a rock.}
\label{f:ar1}
\end{figure}
\v

\v\begin{remark}
{\rm In practice, the rate of growth of the root will likely depend on 
the hardness of the soil at the tip.   In other words, at time $t$ the total length
will not be $\ell(t)=t$.  Rather, it may increase at a rate
\bel{g2}\dot\ell(t)~=~{1\over 1+h\bigl(\gamma(\ell(t))\bigr)} \,,\eeq
so that 
\bel{g3} \dot P(t)~=~{d\over dt} \gamma(t,\ell(t))
~=~{\bfk(\ell(t))\over 1+h\bigl(\gamma(\ell(t))\bigr)} \,.\eeq
This more general  case can be reduced to the previous setting
 by using a rescaled time
$\tau(t) = \ell(t)$.
In this way, as in \cite{07a9e864592f44a99f73752cdcb6e0a8, MR3809025, MR3650339},
we can assume that at time $\tau$ the curve $\gamma(\tau,\cdot)$
has length $\tau$.   This length thus  increases at unit rate.}
\end{remark}

 \subsection{A flexible root}\label{S33}
\label{s:4}
In the previous model the growing root was completely rigid.   Since the previously constructed portions of the root do not move, the equation of growth reduces to a simple ODE.  

 We now consider an alternative model, where the root is allowed to change its shape, in response to 
 obstacles in  the surrounding environment. The instantaneous deformation of the root is determined as a minimizer of an elastic
 deformation energy, plus a cost for displacing the nearby soil.
 
Let $s\mapsto \gamma(t,s)$, $s\in [0,t]$, be an arc-length 
 parameterization of the root at time $t$.
 Moreover,   let $h=h(x)\geq 0$ be  a function
 describing the hardness of the soil at the point $x$, and fix a constant $\alpha>0$.
In practice, one may think to acquire information about the rigidity of the soil in the surrounding of the
 robotic root by means of sensor located along the body of the root.

If the soil  around the tip is very hard, or if 
obstacles are encountered, the root may be forced to bend.
As in \cite{MR3809025, MR3650339}, introducing a field of angular velocities $\omega(t,\cdot)$, the evolution of the curve $\gamma$ is described by
\bel{21}\gamma_t
(t,s)~=~\int_0^s  \omega(t,\sigma)\times\bigl(\gamma
(t,s)-\gamma
(t,\sigma)\bigr)
\, d\sigma,\eeq
while, with the same setting as in~\eqref{tv}-\eqref{Pk}, we have
\bel{22}
\dot P(t)~=~\bfk(t)+ \int_0^t  \omega(t,\sigma)\times\bigl(\gamma
(t,s)-\gamma
(t,\sigma)\bigr)
\, d\sigma.\eeq
As in (\ref{g4}), this equation must be supplemented with a 
boundary condition describing how the orientation of the tip changes in time:
\bel{g44} \dot \bfk(t)~=~\left(\bfu(t) +\int_0^t  \omega(t,\sigma)\, d\sigma
\right) \times \bfk(t) .\eeq
Here the feedback control $\bfu(t)$ can be chosen  as in (\ref{uopt1}), 
in order to approach the target as fast as possible.  Notice that the integral term 
in (\ref{g44}) accounts for the change in the orientation of the tip coming from 
the elastic deformation.

At each given time $t$, following \cite{MR3809025}
we assume that the field of angular velocities $\omega(\cdot)$ 
in (\ref{21})--(\ref{g44}) is determined by 
an instantaneous minimization problem, which we now describe.
\v

\begin{figure}[ht]
\centerline{\hbox{\includegraphics[width=7cm]{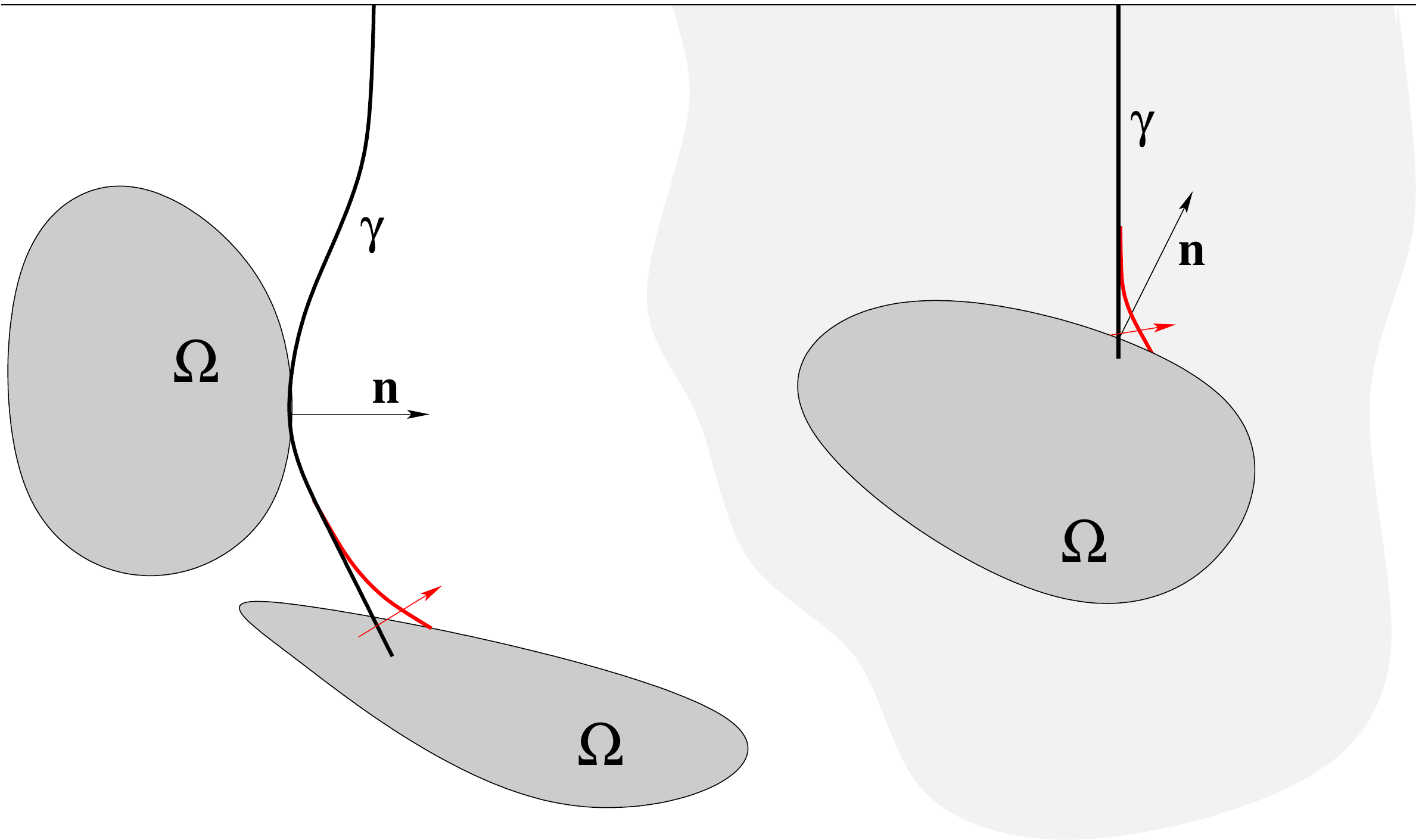}}}
\caption{\small   By an infinitesimal elastic deformation, the curve $\gamma$ is kept
outside the obstacle.  Left: outside the obstacle there is empty space.  Right:
the obstacle is surrounded by softer soil.}
\label{f:ar2}
\end{figure}

Assume that the external environment consists of soil whose hardness
is described by a scalar function $h=h(x)\geq 0$, together
with impenetrable obstacles, whose union we denote by $\Omega\subset\R^3$.  
We assume that $\Omega$ is an open set with $\C^2$ boundary $\partial\Omega$.
The unit outer normal vector at a point $x\in\partial \Omega$ will be denoted
by $\bfn(x)$.  At a contact point $\gamma(t,s)\in \partial\Omega$,
this unit outer normal is denoted by 
$$\bfn(t,s)~\doteq~\bfn(\gamma(t,s)).$$

In our model, the angular velocity 
$\omega(t,\cdot)\in \L^2([0,t])$ is determined as the solution to the following optimization problem.

\begi
\item[{\bf (OP)}] 
{\it Given $\gamma(t,\cdot) \in {\bf H}^2([0,t])$, find an angular velocity $\omega\in \L^2([0,t])$ which minimizes the
functional
\bel{J}\bega{rl}
J(\omega)&\ds\doteq~\int_0^{t} |\omega(s)|^2\, ds + \int_0^t 
 h(\gamma(t,s))\cdot \Big| \gamma_t(t,s)\times \bfk(t,s)\Big|\, ds +\alpha\, h(P(t))\cdot
 \la \dot P(t),\bfk(t)\ra_+\\[4mm]
 &=~J_1(\omega)+J_2(\omega)+J_3(\omega)\,.\enda\eeq
 Here $\gamma_t$ and $\dot P$ are recovered from $\omega(\cdot)$
 by the formulas (\ref{21})-(\ref{22}), while 
 $$\langle \bfv,\bfw\rangle_+~\doteq~\max\{ \langle \bfv,\bfw\rangle\,,~0\}$$
 denotes the positive part of an inner product.
The global minimizer of the functional (\ref{J}) is sought under the constraints
\bel{co1} \Big\langle  \bfn(t,s)\,,\, \gamma_t(t,s)\Big\rangle ~\geq~0\qquad
\hbox{for all $s$ such that $\gamma(t,s)\in \partial\Omega$}.
\eeq
In the case where $P(t)\doteq\gamma(t,t)\in \partial\Omega$, we impose the additional 
constraint
\bel{co2} \Big\langle  \bfn(t,t)\,,\, \dot P(t)\Big\rangle ~\geq~0.\eeq
}
\endi

\begin{remark} \label{rm}{\rm 
The three terms on the right hand side of (\ref{J}) can be interpreted as
$$\hbox{[elastic  bending energy] ~+~ [soil hardness]$\times$[swept area] ~+~
[soil hardness]$\times$[tip penetration].}$$
The coefficient $\alpha>0$ allows to better calibrate the relative weight of these 
terms.}
\end{remark}
\v
\begin{remark} {\rm   In \cite{MR3809025, MR3650339}, the surface of the obstacle was modeled as a
smooth frictionless surface. In alternative, one could impose that there is some friction
between the obstacle and the tip of the root.   For example, one can impose that 
the tip can move only if 
\bel{CF} \la \bfn(t,t)\,,\, \bfk(t,t)\ra~\geq~\kappa_0\,,\eeq
for some $\kappa_0 \in \,] -1,0]$.}
\end{remark}
\v
\begin{remark}  {\rm In the first model, considered in Section~\ref{s:1}, the root is completely rigid, in the sense 
that the portion that is already grown does not change its position at any later time.
On the other hand, in this second model the dynamics is more complex, because it takes into account the elasticity of the root.   

The difference is only in the dynamics.   The optimal control (orienting the tip 
in the direction where it is more likely to find the target), as well as the 
restarting strategy, can be chosen in the same way as in the model discussed 
in the first two sections.}
\end{remark}

 \section{Solutions to the flexible root model}\label{S34}
\label{s:5}
\begin{figure}[ht]
\centerline{\hbox{\includegraphics[width=10cm]{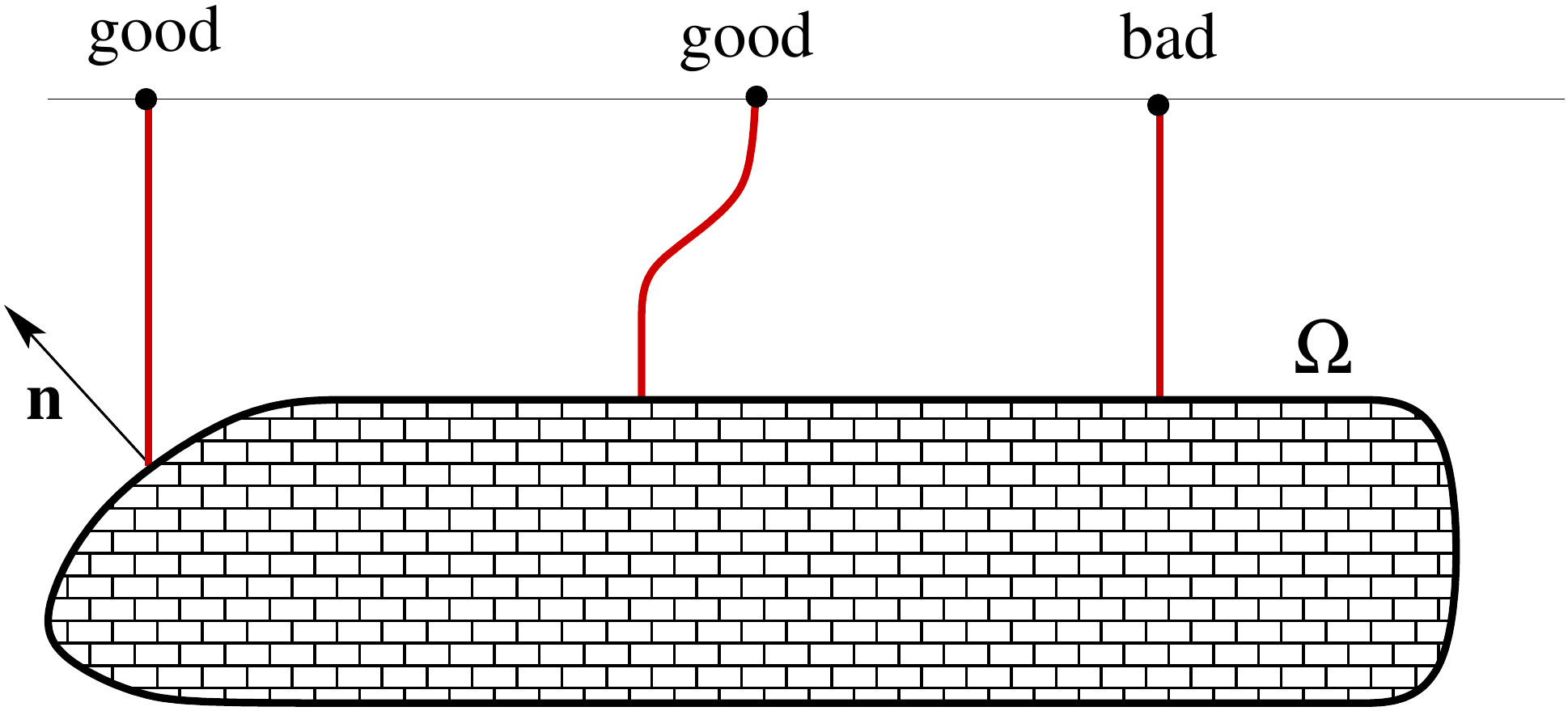}}}
\caption{\small For the two initial configurations on the 
left and in the center, the constrained growth equations for the root admit a unique solution. 
The root configuration on the right satisfies the condition
{\bf (B)}. In such case, the Cauchy problem is ill posed. }
\label{f:sg206}
\end{figure}

\setcounter{equation}{0}
From a theoretical point of view, the model with a flexible root is the most interesting.
Toward a proof of well-posedness, two steps are required.
\v
\begin{itemize}
\item[{\bf (i)}] At any time $t\geq 0$, the instantaneous growth is determined by solving a
 minimization problem with constraints.
One needs to prove that this optimization problem admits a unique solution.
Moreover, this minimizer depends continuously on the data, except 
when a ``breakdown configuration" is reached, as in the model
studied in \cite{MR3650339}.
\v
\item[{\bf (ii)}]  Following \cite{MR3809025}, one should then introduce a suitable weighted distance among root configurations, and
prove that the evolution problem is well posed.
\end{itemize}
\v
Aim of this section is to prove a local existence theorem, valid as long as a 
``breakdown configuration" is not reached, as shown in Fig.~\ref{f:sg206}.

\begin{definition}\label{d:breakdown} We say that a curve $\gamma:[0,t_0]\mapsto \R^3\setminus\Omega$
is in a {\bf breakdown configuration} w.r.t.~the obstacle $\Omega$ if the following holds
\begi
\item[{\bf (B)}]  The tip $\gamma(t_0)$  touches the obstacle perpendicularly, namely
\bel{bad1} \gamma(t_0)~\in~ \partial \Omega\,,\qquad\qquad 
 \gamma_s(t_0)~=~-\bfn(\gamma(t_0)).\eeq 
Moreover, 
\bel{bad2} \gamma_{ss}(s)~= 0\qquad  
\hbox{for all $ s\in ~]0,t_0[\,$ such that~}\gamma(s)\notin \partial\Omega\,.\eeq
\endi
\end{definition}
\v

\v
\subsection{Determining the instantaneous velocity.}

At each time $t\geq 0$, the velocity of points on the curve $\gamma(t,\cdot)$
and of the root tip $P(t)$ 
is determined by (\ref{21}) and (\ref{22}), respectively.    In our model, the instantaneous bending rate $\omega(t,\cdot)\in \L^2([0,t]\,;\,\R^3) $ is determined as the unique global minimizer
of the  functional (\ref{J}), subject to the constraints (\ref{co1})-(\ref{co2}).
We recall that $h:\R^3_-\mapsto\R$ is a smooth scalar function, 
describing the  hardness of the soil.
\v
\begin{lemma}\label{l:51} If $\gamma(t,\cdot)$ is not a breakdown configuration,
then the functional  (\ref{J}) has a unique minimizer, subject to 
(\ref{21})-(\ref{22}) and the constraints (\ref{co1})-(\ref{co2}).
\end{lemma}
\v
{\bf Proof.}  
 {\bf 1.} Call $\A\subseteq\L^2([0,t])$ the set of 
angular velocities for which the constraints (\ref{co1})-(\ref{co2})
are satisfied.    Recalling~\eqref{21}-\eqref{22}, we observe that the maps
$$\omega(\cdot)~\mapsto~\gamma_t(t,\cdot), \qquad\quad \omega(\cdot)~\mapsto~
\dot P(t)$$
are linear and affine, respectively.
We thus deduce that
$\A$ is  a closed, affine subspace of $\L^2$. 

Next,  if $\gamma(t,\cdot)$
does not satisfy all conditions in {\bf (B)},
by adapting some arguments in  the proof of Lemma~1
in~\cite{MR3650339} we now show that $\A$ is  a nonempty set.
Indeed, assume that $\bfv \in {\bf H}^2([0,t];\R^3)$ is a velocity field produced by the obstacle reaction when  $\gamma(t,\cdot)$  touches $\partial \Omega$,  or a reaction due to hardness of the soil, with
\bel{reaction}
	\bfv (0)=\bfv '(0)=0, \quad\quad \langle \bfv'(s)~,~\bfk(t,s) \rangle=0 \quad \text{ for all } s\in[0,t],
\eeq
\bel{reaction2}
\langle \bfn(t,s)\,,\,\bfv(s)\rangle \geq 0
\qquad  \text{if} \qquad
\gamma(t,s)\in \partial\Omega\,,
\quad  s\in[0,t],
\eeq
\bel{reaction3}
\langle \bfn(t,t)\,,\,\bfv(t)+\bfk(t,t)\rangle \geq 0
\qquad  \text{if} \quad\  \gamma(t,t)\in \partial\Omega\,.
\eeq
Here $\bfv$ plays the role of $\gamma_t(t,\cdot)$, defined by~\eqref{21}.
We claim that there exists a unique angular velocity field 
$\omega\in \L^2([0,t];\R^3)$ such that 
\bel{th}
	\Vert \omega \Vert_{\L^2}\leq C \Vert \bfv \Vert_{{\bf H}^2}\,,\\
\eeq
\bel{th0}
	\langle \omega(s), \bfk (t,s) \rangle=0 \quad \quad \text{ for all } s\in[0,t]\,,
\eeq
\bel{th1}
	\bfv(s)= \int_0^s \omega(\sigma)\times (\gamma(t,s)-\gamma(t,\sigma))d\sigma \quad \text{for all } s\in [0,t].
\eeq

\noindent
Notice that \eqref{reaction2}, \eqref{th1} together yield~\eqref{co1}, while
from~\eqref{22}, \eqref{co1} and\eqref{reaction3}, we recover \eqref{co2}.

To prove the Claim, using (\ref{reaction}) and differentiating (\ref{th1}), we first observe that $ \omega(\cdot)$ satisfies  (\ref{th1}) if and only if 
\bel{th2}
	\bfv '(s)~=~ \int_0^s \omega(\sigma)d\sigma \times \bfk(t,s) \quad \text{for all } s\in [0,t].
\eeq

Now  consider a family of orthonormal frames $\lbrace  \bfe_1(s), \bfe_2(s),\bfe_3(s) \rbrace$ such that $ \bfe_1(s)= \bfk(t,s) $ for all $s\in[0,t]$. 
We seek two scalar functions $\omega_2,\omega_3:[0,t]\rightarrow \R$ 
so that  (\ref{th1}) holds with
 \bel{dec}
	\omega(s)~=~\omega_2(s)\bfe_2(s)+\omega_3(s)\bfe_3(s).
\eeq
Notice that (\ref{th0}) is an immediate consequence of the orthogonality of the three vectors 
$\bfe_i(s)$.

By the orthogonality assumption in (\ref{th0}),   
we get two scalar functions $z_2$ and $z_3$ such that 
 \begin{align}
	\bfv'(s)~=~z_2(s)\bfe_2(s)+z_3(s)\bfe_3(s)&=~ \int_0^s \bigl(\omega_2(\sigma)\bfe_2(\sigma)+\omega_3(\sigma)\bfe_3(\sigma)  \bigr)d\sigma \times \bfe_1(s).\label{Dec}
 \end{align}
 Projecting respectively along $\bfe_2(s)$ and  $\bfe_3(s)$ we obtain
 \begin{align*}
	& z_2(s)= \int_0^s \langle \bfe_2(\sigma)  \times \bfe_1(s),~  \bfe_2(s)\rangle \omega_2(\sigma) d\sigma
	+ \int_0^s \langle \bfe_3(\sigma)  	\times \bfe_1(s),~  \bfe_2(s)\rangle \omega_3(\sigma) d\sigma,\\
 	& z_3(s)= \int_0^s \langle \bfe_2(\sigma)  \times \bfe_1(s),~  \bfe_3(s)\rangle \omega_2(\sigma) d\sigma 
	 + \int_0^s \langle 		\bfe_3(\sigma)  	\times \bfe_1(s),~  \bfe_3(s)\rangle \omega_3(\sigma) d\sigma. 
 \end{align*}
By  the property of the mixed product $\langle {\bf u}, ({\bf v} \times {\bf w})\rangle=
\langle {\bf v}, ({\bf w} \times {\bf u})\rangle$,
we can write 
\begin{align*}
	& z_2(s)~=~\ds \int_0^s  \Big[ \langle \bfe_2(\sigma),~ \bfe_3(s)\rangle \omega_2(\sigma) + \langle \bfe_3(\sigma),~ \bfe_3(s)\rangle \omega_3(\sigma)   \Big] d\sigma\,,\\
 	&  z_3(s)~=~\ds -\int_0^s  \Big[ \langle \bfe_2(\sigma),~ \bfe_2(s)\rangle \omega_2(\sigma) + \langle \bfe_3(\sigma),~ \bfe_2(s)\rangle \omega_3(\sigma)   \Big] d\sigma\,.
 \end{align*}
Since all the function involved in $(\ref{Dec})$ are in ${\bf H}^1$, we can differentiate once again  and obtain
\bel{Vie} 
\left\{\begin{array}{rl}
	\omega_3(s)&= ~\ds z_2'(s)-\int_0^s \Big[\langle \bfe_2(\sigma),~ \bfe_3'(s)\rangle \omega_2(\sigma) + \langle \bfe_3(\sigma),~ \bfe_3'(s)\rangle \omega_3(\sigma)   \Big] d\sigma,\\[4mm]
	\omega_2(s)&\ds=~-z_3'(s) -\int_0^s  \Big[ \langle \bfe_2(\sigma),~ \bfe_2'(s)\rangle \omega_2(\sigma) + \langle \bfe_3(\sigma),~ \bfe_2'(s)\rangle \omega_3(\sigma)   \Big] d\sigma.
\end{array}\right.
\eeq
This  is a linear system of of Volterra integral equation, whose unique solution can be obtained by a fix point argument. Indeed, introducing the vector notation  
\bel{COV}
	U=\begin{pmatrix}
		\omega_2\\
		\omega_3
	\end{pmatrix}, \quad \quad 
	Z=\begin{pmatrix}
		-z_3 \\
		z_2
	\end{pmatrix},
\eeq
the system (\ref{Vie}) can be written as 
\bel{CF}
	U(s)~=~Z'(s)+\int_0^sB(s,\sigma)U(\sigma)d\sigma~\doteq
	~ \mathcal{F}[Z',U](s),
\eeq 
where the matrix $B(s,\sigma)$ has norm

\bel{N}
	|B(s,\sigma)|~\leq ~2 |\bfe_2'(s)|+ 2 |\bfe_3'(s)|~\doteq ~b(s).
\eeq
The operator $\mathcal{F}$ defined at (\ref{CF}) is a strict contraction w.r.t.~$U$
on the space $\L^1([0,T]; \mathbb{R}^2)$,
with the equivalent norm
\begin{equation}
	\Vert  U  \Vert~\doteq~
	\int_{0}^t \text{exp}\left\lbrace -2\int_0^s b(\sigma)d\sigma |U(s)|  \right\rbrace |U(s)| ds.
\end{equation}
 Moreover,  by the contraction mapping principle we have
\bel{CH}
\Vert  U  \Vert_{\L ^1} \leq C_1 \Vert  Z' \Vert_{\L ^1}\leq C_2 \Vert   Z'  \Vert_{\L ^2},
\eeq
with $C_1,C_2$ constants depending on $t$ and $b(\cdot)$. 
Conclusively we obtain that 
\begin{align*}
	\Vert  U  \Vert_{\L ^2}^2 &\leq  2\Vert   Z'  \Vert_{\L ^2}^2 + \int_0^t b^2(s)\Big( \int_0^sU(\sigma) \sigma \Big)^2ds\\
	& \leq 2\Vert  Z' \Vert_{\L ^2}^2 + 2\Vert  b\Vert_{\L ^2}^2\Vert  U  \Vert_{\L ^1} ^2\\
	&\leq (2+ 2\Vert  b\Vert_{\L ^2}^2C_2^2)\Vert  Z'  \Vert_{\L ^2}\\
\end{align*}
hence also $\omega_1$ and $\omega_2$ are in $\L^2$.  To conclude we just need to observe that $\Vert   Z'  \Vert_{\L ^2}=\mathcal{O}(1)\cdot \Vert \bfv \Vert_{{\bf H}^2}$ and $\Vert b \Vert_{\L^2} \mathcal{O}(1)\Vert \gamma\Vert_{{\bf H}^2} $.
Condition $(\ref{th0})$ plays an essential role in uniquess of the angular velocity $\omega$, indeed without that we can produce infinitely many solution for the system $(\ref{th1})$.

\v
{\bf 2.} We now observe that the functional $J$ in (\ref{J}) is non-negative, strictly convex, and coercive on $\L^2$.   

Let  $(\omega_n)_{n\geq 1} $ be a minimizing sequence in $\A\subseteq\L^2$.
By coercivity, the norms $\|\omega_n\|_{\L^2}$ are uniformly bounded.
By possibly extracting a subsequence, we can thus assume the weak convergence
$\omega_n\wto\omega^*$.  

Since the functional $J$ is convex, it is  lower semicontinuous w.r.t.~weak convergence. Therefore
$$J(\omega^*) ~\leq~\liminf_{n\to\infty} J(\omega_n)~=~\inf_{\omega\in \A} J(\omega).$$
 Recalling that $\A$ is closed, we conclude that $\omega^*\in\A$ is a global minimizer.
 
 Finally, the strict  convexity of $J$ implies that this minimizer is unique.
 \endproof
 
 \v

 \subsection{Solutions to the evolution problem.}
 
 Given a feedback control $\bfu=\bfu(x)$ accounting for the desired bending of the tip, 
we can now give a precise definition of ``solution" of the evolution problem 
 describing the growth of the artificial root.

At each time $t$, the position of the root is described by a map
$\gamma(t,\cdot)$ from $[0,t]$ into $\R^3$.    Of course, the domain of this map grows with
time.   Following  \cite{MR3650339},
it is convenient to reformulate our model as an evolution problem 
on a functional space independent of $t$.   We thus fix $T>t_0$
and consider the Hilbert-Sobolev space ${\bf H}^2([0,T];\, \R^3)$.  

An initial data $\overline\gamma\doteq \gamma(t_0, \cdot) \in {\bf H}^2([0,t_0];\, \R^3)$ 
will be canonically extended to 
${\bf H}^2([0,T];\, \R^3)$ by setting 
\bel{Pex}
\overline\gamma(s)~\doteq~\overline\gamma(t_0)+(s-t_0)\overline\gamma'(t)
\qquad\qquad 
\hbox{for}~~s\in [t_0,T]\,.\eeq
Notice that the above extension is well defined because 
$\ov \gamma(\cdot)$ and its derivative $\ov \gamma'(\cdot)$
are continuous functions. In the following we
shall study
functions defined on a domain of the form
\bel{DT}\D_T~\doteq~\bigl\{(t,s)\,;~~0\leq s\leq t,~~
t\in [t_0,T]\bigr\}, \eeq
and extended to the rectangle $[t_0,T]\times [0,T]$ as in (\ref{Pex}).

 
\begin{definition}\label{d:root}
We say that a function $\gamma=\gamma(t,s)$, defined for 
$(t,s)\in [t_0, T]\times [0,T]$ is a solution to the root growth problem~\eqref{21}-\eqref{g44},
with initial condition
\bel{id} \gamma(t_0, s)~=~\ov \gamma(s),\qquad\qquad s\in [0,t_0],\eeq
 if the following holds.
\begi
\item[(i)] The map
$t\mapsto P(t,\cdot)$ is Lipschitz continuous from $[t_0,T]$ into ${\bf H}^2([0,T];\,\R^3)$.
\item[(ii)]  For a.e.~ $t\in [t_0,T]$ the partial derivative
$\partial_t \gamma(t,\cdot)$ is given by 
\bel{gt1}\gamma_t
(t,s)~=~\int_0^s  \omega(t,\sigma)\times\bigl(\gamma
(t,s)-\gamma
(t,\sigma)\bigr)
\, d\sigma,\eeq
for $s\in [0,t]$, while
\bel{gt2}\gamma_t
(t,s)~=~\int_0^t  \omega(t,\sigma)\times\bigl(\gamma
(t,s)-\gamma
(t,\sigma)\bigr)
\, d\sigma~+~\bfu(t) \times\bigl(\gamma
(t,s)-\gamma
(t,t)\bigr)\eeq
if $s\in [t,T]$.
Here the angular velocity $\omega(t,\cdot)$
is the unique solution to the optimization problem {\bf (OP)}, in Section~\ref{s:4}.
 \item[(iii)] The initial conditions hold:
\bel{ic4}
\gamma(t_0,s)~=~\left\{\bega{cl}\ov \gamma(s)&\qquad\hbox{if}~~s\in [0, t_0],\\[4mm]
\ov \gamma(t_0)+(s-t_0)\ov \gamma'(t_0)&\qquad\hbox{if}~~s\in [t_0,T].
\enda\right.\eeq
\item[(iv)]
The pointwise constraints hold:
\bel{pc4}
 \gamma(t,s)~\notin~\Omega\qquad\quad\forall t\in [t_0,T],~~s\in [0,t].\eeq
 \endi
\end{definition}

\section{Construction of solutions}\label{S35}
\subsection{Technical lemmas and known estimates}
We now recall some notations, technical results and estimates contained in  \cite{MR3650339} that will be useful in order to prove the main theorem of the paper. \\

\smallskip
Let $\ov{\gamma}$ be a curve in ${\bf H}^2([0,t_0];\R^3)$, satisfying the following assumptions

\bel{Cent}
	\ov{\gamma}(0)=0\in\R^3,\quad \quad \vert\ov{\gamma}'(s)\vert=1,\quad\quad \ov{\gamma}(s)\notin \Omega \quad \forall s\in[0,t_0].
\eeq

Moreover, assume that $\ov\gamma$ does not satisfy the breakdown conditions {\bf (B)} in Definition~\ref{d:breakdown}. Given $T>t_0$, we 
extend  $\ov{\gamma}$ to a map $[0,T]\rightarrow \R^3$ as in
(\ref{Pex}).
For a fixed radius $\rho$ and $T>t_0$ we can define the following tubular neighborhood of $\ov{\gamma}$ in ${\bf H}^2([0,T],\R^3)$  
\bel{tub}\bega{rl}
	\mathcal{V}_{\rho} &\ds\doteq~ \Big\lbrace \gamma \in ~{\bf H}^2([0,T],\R^3), \quad \gamma(0)=0, \quad \gamma'(0)=\ov{\gamma}'(0),\\[4mm]
	&\qquad\quad \ds\vert \gamma'(s) \vert =1  \text{ for all } s\in[0,T], \quad \int_0^T \bigr| \gamma''(s)-\ov{\gamma}''(s)\bigl|^2 ds\leq \rho \Big\rbrace.
\enda
\eeq
The following result shows that, given a curve satisfying $(\ref{Cent})$ and which it is not  in a breakdown configuration, there exists a tubular neighborhood as in
(\ref{tub}) such that each curve $\gamma\in \mathcal{V}_{\rho}$ can be pushed away from the obstacle  by a bounded deformation.

\begin{lemma}
\label{le}
Let $\ov \gamma:[0,t_0]\mapsto \R^3\setminus\Omega$ an initial curve which satisfies $(\ref{Cent})$ and  does not satisfies simultaneously all conditions in (\ref{bad1})-(\ref{bad2}). Then there exists $T>t_0$, and constants $C_0, \rho,\delta>0$ such that the following holds.\\

Let $\phi :\R^3\mapsto \R^+$ 
be the  signed distance function from $\partial\Omega$, defined as 
\bel{Sd}
	\phi(x)~\doteq~ \begin{cases}
	d(x,\partial\Omega) &\text{ if } x\notin\Omega,\\
	-d(x,\partial\Omega) &\text{ if } x\in\Omega.
	\end{cases}
\eeq
Let $t\in [t_0, T]$ and consider any curve $\gamma\in\mathcal{V}_\rho$. 
Then there exists $\omega:[0,t]\mapsto \R^3$, with 
\bel{Bo}
	\Vert \omega \Vert_{\L^2}\leq C_0
\eeq
such that, for all $s\in[0,t]$ with $\vert \phi(\gamma(s)) \vert \leq \delta $ one has 
\bel{pus}
	\Big\langle \int_0^s \omega(\sigma)\times (\gamma(s)-\gamma(\sigma)),\nabla \phi (\gamma(s)) \Big\rangle \geq 1\,.
\eeq
\end{lemma}
\noindent
For a proof of Lemma~\ref{le} see  \cite[Lemma~2]{MR3650339}. 
\v
In connection with any  curve $\gamma: [0,T]\to \R^3$, which possibly intersects $\Omega$,
we  now introduce the quantity
\bel{depth}
	E(t,\gamma, \Omega)~\doteq~ \sup \Big\lbrace d(\gamma(s),\partial \Omega); \quad s\in[0,t],\gamma(s)\in\Omega\Big \rbrace .
\eeq

\noindent
This quantity measures the maximum depth at which the first portion of the curve $\gamma$ of lenght~$t$
penetrates inside the obstacle $\Omega$.
Since the boundary $\partial\Omega$ is smooth,
in the same hypothesis of Lemma~\ref{le} one can show that there exists a tubular neighborhood $\mathcal{V}_{\rho}$ of $\ov{\gamma}$ such that, for every curve $\gamma$ inside $\mathcal{V}_{\rho}$, 
there exists a unique angular velocity $\ov{\omega}$ minimizer  of $\textbf{OP}$ which pushes out the curve (i.e. such that~\eqref{pus} holds),
and its norm is bounded by 
\bel{bo1}
	\Vert \ov{\omega} \Vert_{\L^2}~\leq~ 2C_0 \cdot E(t,\gamma, \Omega),
\eeq
for some constant $C_0$ independent of $\gamma \in\mathcal{V}_\rho$
(see  \cite[Lemma~3]{MR3650339}).

Next, we introduce a nonlinear ``push-out" operator $\mathcal{P}$.
 Given an angular velocity  $\omega\in\R^3$, let $R[\omega]$ be the $3\times3$ rotation matrix 

\bel{rot}
	R[\omega]~\doteq~e^{A}~=~\sum_{k_0}^\infty \frac{A^k}{k!}, \quad \quad
	 A\,\doteq\,\begin{pmatrix}
0 & -\omega_3 & \omega_2 \\
\omega_3 & 0 & -\omega_1 \\
-\omega_2 & \omega_1 & 0
\end{pmatrix}.
\eeq
Notice that, for every $\ov\bfv\in\R^3$, 
one has 
\bel{sys2}
A \ov\bfv = \omega\times \bfv\,.
\eeq
If $\gamma$ is a curve in ${\bf H}^2([0,t],\R^3 )$, and $\ov{\omega}$ is the unique minimizer of the associated problem $\textbf{(OP)}$, we then define 
\bel{nlr}
	\mathcal{P}^{\,\ov{\omega}}[\gamma](s)~\doteq~
	\int_0^s R\Big[\int_0^\sigma \ov{\omega}(\xi)d\xi \Big]\gamma'(\sigma) d\sigma.
\eeq  

\noindent
By $(\ref{bo1})$, the depth at which $\mathcal{P}^{\,\ov{\omega}}[\gamma]$ can still 
penetrate inside $\Omega$  is bounded by
\bel{depth1}
E(t,\mathcal{P}^{\,\ov{\omega}}[\gamma], \Omega)\leq  \mathcal{O}(1) \Vert \ov{\omega}\Vert_{\L^2}^2\leq C \cdot E^2(t,\gamma, \Omega),
\eeq
for some  uniform constant $C$.

\subsection{An existence result.}
Here we state and prove the main result of the paper. Namely, the existence of a solution to the evolution problem for the root, up to the first time
where a breakdown configuration is reached. 
We remark that, in practice,  the growth of an artificial root
can be interrupted also because of the physical limits of the machine itself, like a bound on its total length.
\begin{theorem}\label{t:exist} Let 
$u= u(P)$ be a Lipschitz continuous feedback control function.
Assume that the initial configuration $\ov\gamma\in {\bf H}^2([0, t_0];\R^3)$
is not a breakdown configuration according with Definition~\ref{d:breakdown}. 
Then there exists $T>t_0$ and a solution
to the root evolution problem (\ref{21})-(\ref{g44}), with initial condition $\overline\gamma$,
defined on an interval $[t_0,T]$.
Moreover, either $T=+\infty$, or else, as $t\to T-$, a breakdown configuration
satisfying all conditions in {\bf (B)} is reached.
\end{theorem}

{\bf Proof.} {\bf 1.}
Differentiating (\ref{gt1})-(\ref{gt2}) w.r.t.~$s$ and then integrating w.r.t.~$t$,  recalling~\eqref{tv} we obtain
\bel{eqt}
	\bfk(t,s)~=~ \bfk(t_0,s)+ \int_{t_0}^t \Big( \int_0^{\tau\wedge s} \omega(t,\sigma)d\sigma \Big) \times \bfk(\tau,s)d\tau + \mathcal{H}(s-t_0) \int_{t_0}			^{t\wedge s} u\big(\gamma(\tau,\tau)\big)\times \bfk(\tau,s)d\tau,
\eeq
with  $(t,s)\in [t_0,T]\times [0,T]$, where $\mathcal{H}$ denotes the Heaviside function and we adopt the notation
$a\wedge b = \min\{a,b\}$.
\v
{\bf 2.} 
A sequence of approximate solutions will be constructed using an operator 
splitting scheme, and relying on the integral equation~\eqref{eqt} describing the evolution of the
tangent vectors to the solution of~\eqref{21}-\eqref{g44}.
Fix a time step $\varepsilon>0$ and set $t_k=t_0+k\varepsilon$, $k\in\N$. Assume that the 
curves $\gamma=\gamma(t,s)$ have been constructed for
all times $t\in[0,t_k]$ and for $s\in[0,t]$. Call 
$\bfk~=~\bfk (t,s)= \gamma_s(t,s)$ the corresponding tangent vectors. 
As in (\ref{ic4}), we extend each curve $\gamma(t,\cdot)$ by a straight 
line for $s\in (t,T]$.  The tangent vectors thus satisfy
\bel{vi}
		\bfk(t,s) ~=~\bfk(t, t)\qquad\hbox{for}\quad  s\in[t, T].
\eeq

To extend the solution on $[t_k,t_{k+1}]$ we proceed in two stages.

\noindent
\textsc{Stage 1:} {\it growth of the tip, whose bending is determined by the
feedback function $u$.}

\bel{SO}
	\bfk (t,s) ~=~  R\left[ ((t\wedge s)-t_k)u\big(\gamma(t_k,t_k)\big)\mathcal{H}(s-t_k) \right] \bfk(t_k,s)
\eeq
 for all $t\in[t_k,t_{k+1} [$ and $s\in [0,T]$, with 
 \bel{previouscurve}
 \nonumber
 \gamma(t_{k},t_k) =\int_0^{t_k} \bfk (t_k,\sigma) d\sigma\,.
 \eeq
 Let $\ov\omega_k\in\L^2([0,t_k])$ be the angular velocity obtained as solution of $(\bf{OP})$ at time $t=t_k$,
 in connection with the curve $\gamma(t_k,\cdot)\in{\bf H}^2([0,t_k])$. Taking $t=t_{k+1}$, the previous construction produces a curve 
 
 \bel{curve}
 	s~ \mapsto~ \gamma(t_{k+1}-,s) =\int_0^s \bfk (t,\sigma) d\sigma.
 \eeq

\noindent
\textsc{Stage 2:} {\it push-out from the obstacle.}

\noindent
Observe that the physical meaningful portion of the curve, that is the part with $s\in[0,t_{k+1}]$ may lie inside the obstacle therefore we need to use a push-out operator and replace $\gamma(t_{k+1}-,s)$ by a new curve, by setting

\bel{push-out k+1}
	\gamma (t_{k+1} , s)~=~ \mathcal{P}\left[  \gamma(t_{k+1}-,\cdot) \right](s)~=~ \int_0^s 					R\left[  \int_0^{\sigma}\ov\omega_{k+1}(\xi)d\xi \right]   \gamma_s (t_{k+1} -, \sigma) d\sigma
\eeq 
\noindent
with  $\ov\omega_{k+1}\in\L^2([0,t_{k+1}])$  solution to $\bf (OP)$ at time $t=t_{k+1}$,
in connection with the curve $\gamma(t_{k+1},\cdot)\in{\bf H}^2([0,t_{k+1}])$. This is equivalent to say that 

\bel{pt}
\begin{aligned}
	\bfk (t_{k+1} , s)~&=~ 
	R\left[  \int_0^{s\wedge t_{k+1}}\ov\omega_{k+1}(\xi)d\xi \right]   \bfk (t_{k+1}- , s)
	\\
	&=R\left[  \int_0^{s\wedge t_{k+1}}\ov\omega_{k+1}(\xi)d\xi \right] \circ
	R\left[ ((t_{k+1}\wedge s)-t_k)u\big(\gamma(t_k,t_k)\big)\mathcal{H}(s-t_k) \right] \bfk(t_k,s)
	\end{aligned}
	\eeq
\noindent
for all $s\in [0,T]$.
By $(\ref{bo1})$ and $(\ref{depth1})$, as long as the approximation remains inside the neighborhood $\mathcal{V}_\rho$, there exists a constant  $C_3$ such that 
\bel{ex1}
	\Vert \ov{\omega}_{k+1} \Vert_{\L^2([0,t_k+1])}\leq C_3 E\big(t_{k+1}, \gamma(t_{k+1}-, \cdot),\Omega\big),
\eeq
\bel{ex2}
	E\big(t_{k+1}, \gamma(t_{k+1}, \cdot),\Omega\big)\leq C_3 E^2\big(t_{k+1}, \gamma(t_{k+1}-, \cdot),\Omega\big).
\eeq
This last inequality express the fact that every time we apply the non linear push-out $(\ref{push-out k+1})$, this entails a penetration of the  curve of $\varepsilon ^2$ if the length of the curve inside the obstacle is of the order~$\varepsilon$. 

%
When we apply the rotation $(\ref{SO})$, the elongation of the tip
%
%
increases of a quantity of the order  of $\varepsilon$, that is 
\bel{DI}
	E\big(t_{k+1}, \gamma(t_{k+1}-,\cdot),\Omega\big)~\leq~ E\big(t_{k}, \gamma(t_{k},\cdot),\Omega\big)+ C_4 \varepsilon
\eeq
for some constant $C_4$. 
When $\varepsilon>0$ is small enough, the estimates  $(\ref{ex2})$-$(\ref{DI})$ 
thus yield the implication
\bel{impl}
	E\big(t_{k}, \gamma(t_{k},\cdot),\Omega\big)~\leq ~\varepsilon\quad  \implies  \quad E\big(t_{k+1}, \gamma(t_{k+1}, \cdot),\Omega\big)	~\leq ~\varepsilon.
\eeq
Since by assumption  the initial curve $\ov{\gamma}$ lies outside the obstacle, 
one has $E\big(t_0,\ov{\gamma}(t_{0},\cdot),\Omega\big)=0$.
By induction, for all $k\geq 1$ we thus have 
\bel{depth2}
	E\big(t_{k+1}, \gamma(t_{k+1}, \cdot),\Omega\big)~\leq ~\varepsilon.
\eeq
Combining (\ref{depth2})  with (\ref{ex1})-(\ref{ex2}), one obtains
\bel{sav}
	\Vert \ov{\omega}_{k+1} \Vert_{\L^2 [0,t_{k+1}]} ~\leq ~C_5 \varepsilon.
\eeq
\v
{\bf 3.} By the previous construction, for every time step $\varepsilon>0$ small enough, we obtain an 
approximate solution $\bfk_{\varepsilon}=\bfk_{\varepsilon}(t,s)$ of~\eqref{eqt},
defined for all $s\in[0,T]$ and $t\in[t_0,T_{\varepsilon}]$, where $T_{\varepsilon}$ is given by
\begin{equation*}
	T_{\varepsilon}=\sup \Big\lbrace \tau\in[t_0,T]: ~ \gamma_{\varepsilon}(\tau, \cdot) \doteq \int_0^{(\cdot)} \bfk_{\varepsilon}(\tau,\sigma)d\sigma \in 		\mathcal{V}_\rho 
	\Big \rbrace ,
\end{equation*}
with $\mathcal{V}_\rho $ denoting the tubular set in~\eqref{tub}.

As long as the approximation $\gamma_{\varepsilon}(\tau, \cdot)$ remains inside $\mathcal{V}_\rho$ we have the following estimates:
\bel{ST1}
	\Vert   \bfk_{\varepsilon}(t,\cdot)-  \bfk_{\varepsilon}(t',\cdot)\Vert_{{\bf H}^1([0,t_k])}~\leq ~C_6 |t-t'|\qquad \forall t,t'\in[t_k,t_{k+1}],
\eeq
\bel{ST2}
	\Vert   \bfk_{\varepsilon}(t_{k+1},\cdot)-  \bfk_{\varepsilon}(t_{k+1}-,\cdot)\Vert_{{\bf H}^1([0,t_{k+1}])}~\leq ~C_6\varepsilon,
\eeq
for some constant $C_6$ independent of $k, \varepsilon$. The first inequality 
follows from the boundedness of $u\big(\gamma(t,\cdot)\big)$, 
while the second one is a consequence of  $(\ref{sav})$. 
Combining (\ref{ST1}) and (\ref{ST2}) we obtain by induction that, 
as long as  $\gamma_{\varepsilon}(\tau, \cdot)\in \mathcal{V}_\rho$ we have
\bel{ST3}
	\Vert   \bfk_{\varepsilon}(\tau,\cdot)-  \bfk_{\varepsilon}(t,\cdot)\Vert_{{\bf H}^1([0,t])}\leq C_7 (m\,\varepsilon +\tau-t) ,\eeq
for  some constant $C_7$, and all $t_h\leq t<\tau \leq t_{k+1}$, $m= [(t_{k+1}-t_h)/\varepsilon]$,
where $[\cdot]$ denotes the integer part.


By construction, when $s=0$ we always have
\bel{ST0}
	\gamma_{\varepsilon}(t,0)\,=\,0,
	\qquad \quad  \bfk_{\varepsilon}(t,0)=\ov{\bfk}(0) \quad \forall t\geq 0. 
\eeq
Therefore, (\ref{ST3}) implies 
\bel{ST4}
	\gamma_{\varepsilon}(t,\cdot) \in\mathcal{V}_\rho \quad \forall t\in[t_0,T],
\eeq
for some $T>t_0$ sufficiently close to $t_0$, depending on $\rho$ but independent of $\varepsilon$. The estimates implies also that, for $\varepsilon$ small enough, all the approximations  $\bfk_{\varepsilon}= \bfk_{\varepsilon}(t,s)$  are well defined.
\v

{\bf 4.} To prove 
 the convergence of a subsequence of approximations we observe that:
\begi
\item[(i)]  All the functions $t\mapsto \bfk_{\varepsilon}(t,\cdot)$ have uniformly bounded total variation, as maps from $[t_0,T]$ into ${\bf H}^1([0,T])$.
\item[(ii)]
 Since $\gamma_{\varepsilon}(t,\cdot)\in \mathcal{V}_\rho$,  for any fixed $t\in[t_0,T]$ we have that 
 $\{\bfk_{\varepsilon}(t,\cdot),\, \varepsilon>0\}$, are uniformly bounded in ${\bf H}^1([0,T])$. 
Therefore, by a version of Helly's selection principle for BV functions with values in metric spaces (see \cite{MR2210910, MR2105969}), 
there exists a subsequence  $\bfk_{\varepsilon_n} $ and a function $\bfk \in BV([t_0,T]\,;~ {\bf H}^1([0,T]),\R^3)$ such that  $\bfk_{\varepsilon_n}(t,\cdot)\rightharpoonup \bfk(t,\cdot)$ in ${\bf H}^1([0,T])$ for every $t\in[t_0,T]$. By the compact embedding of ${\bf H}^1([0,T])$ into $C([0,T])$ we can  also achieve the uniform convergence  $\bfk_{\varepsilon_n}(t,s)\to \bfk(t,s)$ on $[t_0,T]\times [0,T]$.
\endi
\v

{\bf 5.}
Aim of this step is to prove that  the limit function $\bfk(\cdot,\cdot)$ satisfies
(\ref{eqt}).
Performing an integration, the convergence $\bfk_{\ve_n}\to\bfk$
implies the  uniform  convergence
\bel{cov}
\gamma_{\varepsilon_n}(t,s)~=~\int_0^s \bfk_{\varepsilon}(t.\sigma)d\sigma ~\rightarrow~ \int_0^s \bfk(t.\sigma)d\sigma ~=~ \gamma(t,s).
\eeq

For every $t\in[t_0,T]$, consider the sequence of convex  functionals
\bel{fcs}
	J_{\varepsilon}(\omega)=~\int_0^{t} |\omega(s)|^2\, ds + \int_0^t 
	 h(\gamma_{\varepsilon}(t,s))\cdot \Big| (\gamma_{\varepsilon})_t(t,s)\times \bfk_{\varepsilon}(t,s)\Big|\, ds +\alpha\, h(P_\varepsilon(t))\cdot
	 \la \dot P_{\varepsilon}(t),\bfk_{\varepsilon}(t,t)\ra_+
\eeq
with $\omega\in\L^2([0,t])$.
 These functionals are equi-bounded in bounded sets of $\L^2([0,t])$.
 
As $\ve_n\to 0$, they converge pointwise to the convex functional $J$ in (\ref{J}). By Theorem 5.12  in \cite{MR1201152}, the sequence  
$( J_{\varepsilon_n})_{n\geq 1}$
 $\Gamma$-converges to $J$ in $\L^2([0,t])$ as $\varepsilon_n \rightarrow 0$. Hence we have also the convergence of the sequence of minimizers $\omega^{\varepsilon_n}(t,\cdot)$ of $J_{\varepsilon_n}$
 to the unique minimizer $\omega(t,\cdot)$ of $J$, in the $\L^2([0,t])$ topology. The $\L^2$-convergence in a bounded domain implies the $\L^1$ convergence and hence the pointwise convergence  of $\omega^{\varepsilon_n}(t,\cdot)$
 a.e. in $[0,t]$.


Using  the estimate
\begin{equation*}\Big\vert \bfv + (\omega_1 + ...+\omega_n)\times \bfv - R[\omega_1]\circ \cdot \cdot \cdot\circ R[\omega_n] \bfv\Big\vert=\mathcal{O}(1) \Big(\sum_i |\omega_i|\Big)^2 |\bfv|,
\end{equation*}
together with the uniform convergence $\bfk_{\varepsilon}\to\bfk$ and  the  pointwise a.e.~convergence $\omega_{\varepsilon}\to \omega$, and recalling~\eqref{eqt}, \eqref{pt}, we conclude that 
$$\bega{rl}
	\bfk(t,s)-\bfk(t_0,s)&\ds= ~\lim_{\varepsilon\rightarrow 0}  \mathcal{H}(s-t_0)\sum_{k=0}^{k_{\varepsilon}(t)} (t_{k+1}\wedge s -t_k)u(t_k)\times 		\bfk_{\varepsilon}(t_k^{\varepsilon},s)\\[4mm]
	&\ds\qquad +\lim_{\varepsilon\rightarrow 0}\sum_{k=0}^{k_{\varepsilon}(t)} \Big(\int_{0}^{t_k\wedge s}\omega^{\varepsilon}_{k}(\sigma)d\sigma\Big) 		\times\bfk_{\varepsilon}(t_k^{\varepsilon},s),
\enda $$
where  $t_k^{\varepsilon}=t_0+k\varepsilon$ and $k_{\varepsilon}(t)= \max \lbrace k\geq0,~t_k^{\varepsilon}\leq t \rbrace$.
\v
{\bf 6.}  To complete the proof we
 show that the curve $\gamma$ defined by
 \bel{gammacurve}
 \gamma(t,s) ~=~\int_0^s \bfk (t,\sigma) d\sigma\,,
 \eeq
is differentiable for a.e. $t\in[t_0,T]$   as a map with values in ${\bf H}^2$.

By the previous steps, and because of the integral representation~\eqref{eqt} of $\bfk$,
the derivative $\bfk_{t}$ is well defined for a.e. $(\tau,s)\in[t_0,T]\times [0,T]$ and it satisfies a uniform bound 
with respect to~$s$ $|\bfk_{t}|\leq C_t$.  
Hence, by the Lebesgue dominated convergence theorem, for a.e. $t\in[t_0,T]$ we obtain
\begin{align*}
\gamma_t(t,s)&= \lim_{\varepsilon\rightarrow 0} \frac{\gamma(t+\varepsilon ,s)-\gamma(t,s)}{\varepsilon}= \lim_{\varepsilon\rightarrow 0}  
\int_0^s \frac{\bfk(t+\varepsilon ,\sigma)-\bfk(t,\sigma)}{\varepsilon}d\sigma= \int_0^s \bfk_t(t,\sigma)d\sigma\,.
\end{align*}
This implies the Lipschitz continuity of $\gamma$ as map from $[t_0,T]$ to ${\bf H}^2([0,T],\R^3)$. 

\v
{\bf 7.} Set 
\bel{t-break-def}
\overline T\doteq \sup\Big\lbrace
t>t_0:~\gamma(t,\cdot)\big|_{[0,t]}~~\text{is not a breakdown configuration}
\Big\rbrace,
\eeq
and consider an increasing sequence $t_n>t_0$, with $t_n\uparrow \overline T$.
By the previous costruction there will be times $T_n\in (t_n, \ov T\,)$, and curves
$\gamma_n: [t_0, T_n]\times [0,T_n]\to \R^3$,  which are
solutions to the root evolution problem (\ref{21})-(\ref{g44}), with initial condition $\overline\gamma$,
defined on the intervals $[t_0,T_n]$.
Notice that $T_n\uparrow \overline T$. Then, 
setting 
$$
\bfk_n(t,s)=(\gamma_n)_s(t,s)\qquad \forall~(t,s)\in [t_0, T_n]\times [0,T_n],
$$
with the same arguments of the above points we  find:
\begi
\item[(i)]  Because of~\eqref{ST3},
all the functions $t\mapsto \bfk_n
(t,\cdot)$ have uniformly bounded total variation, as maps from $[t_0,T_n]$ into ${\bf H}^1([0,T_n])$.
\item[(ii)]
 Since $\bfk_n$ satisfy~\eqref{eqt} it follows that,  for any fixed $t\in[t_0,\overline T)$,
 there holds
 \bel{h1-unif-bound}
 \sup_n \Big\{\big\|\bfk_n(t,\cdot)\big\|_{{\bf H}^1([0,T_n])}: t<T_n\Big\}<\infty\,.
 \eeq
Therefore, relying on the  Helly's selection principle for BV functions with values in metric spaces 
and on the compact embedding of ${\bf H}^1([0,T_n])$ into $C([0,T_n])$ we 
deduce that there is a function
$\bfk \in BV([t_0,\overline T)\,;~ {\bf H}^1([0,\overline T]),\R^3)$ such that
a subsequence of $\bfk_n$ converges to~$\bfk$, uniformly on compact subsets of $[0,\overline T)\times [0,\overline T)$.
\item[(iii)]
The convex  functionals
\bel{fcs}
	J_n(\omega)=\int_0^{T_n} |\omega(s)|^2\, ds + \int_0^t 
	 h(\gamma_n(t,s))\cdot \Big| (\gamma_n)_t(t,s)\times \bfk_n(t,s)\Big|ds +\alpha\, h(P_n(t))\cdot
	 \la \dot P_n(t),\bfk_n(t,t)\ra_+
\eeq
with $\omega\in\L^2([0,\widehat T])$, are equi-bounded in bounded sets of $\L^2([0,\overline T])$.
Therefore, we deduce the $\Gamma$-convergence in $\L^2([0,\ov T])$ of $J_n$ to the functional $J$ in~\eqref{J}.
Hence,  we have also the convergence of a sequence of 
minimizers $\omega_n(t,\cdot)$ of $J_n$
 to a minimizer $\omega(t,\cdot)$ of $J$, in the $\L^2([0,\overline T])$ topology.
\item[(iv)]
Notice that $\omega_n(t,\cdot)\big|_{[0,T_n]}$ are the unique minimizer of $J_n$ in the $\L^2([0,T_n])$
space. Therefore, by~\eqref{eqt} they satisfy
 \bel{eqt2}
 \begin{aligned}
	\bfk_n(t,s)~&=~ \bfk_n(t_0,s)+ \int_{t_0}^t \Big( \int_0^{\tau\wedge s} \omega_n(t,\sigma)d\sigma \Big) \times \bfk_n(\tau,s)d\tau + 
	\\
	& +\mathcal{H}(s-t_0) \int_{t_0}^{t\wedge s} u\big(\gamma_n(\tau,\tau)\big)\times \bfk_n(\tau,s)d\tau,
	\qquad\quad \forall~(t,s)\in [t_0,T_n]\times[0,T_n]\,.
	\end{aligned}
\eeq
then, passing to the limit as $n\to \infty$ in~\eqref{eqt2}, we deduce that $\bfk$ is a solution of~\eqref{eqt},
which implies that the curve defined by~\eqref{gammacurve} provides a solution 
to the root evolution problem (\ref{21})-(\ref{g44}) on $[t_0,\overline T)\times [0,\overline T)$.
\endi

\section{Numerical simulations}\label{S36}
In this final section we present some numerical simulations for the 
flexible root model proposed in Section~\ref{S33}. 
We first show how to rewrite it in the planar case, and the simulations will be done in dimension two. For any vector $\bfv=(v_1,v_2)$ let $\bfk^\perp = (-v_2,v_1)$ be the perpendicular vector obtained by a counterclockwise rotation of $\pi/2$. Then the evolution equation for the tangent vectors to the curve $\gamma$ 
is given by \bel{pte}
\bfk_t(t,s)= \Big( \int_0^{t\wedge s} \omega(\sigma)d\sigma + \mathcal{H}(s-t)u(t)\Big) \bfk^\perp(t,s).
\eeq
\label{s:15}
\setcounter{equation}{0}
The root is discretized with uniform arc-length $\Delta s$ and each time step is taken such that $\Delta t=\Delta s$. Given an initial configuration, represented by an array of nodes, we construct a new node following the direction of the target. In our simulation the target will be the whole axis $y=0$, which means that the first goal of the root  is to  grow in depth compatibly with the prescribed bound on the curvature, that is, the root cannot suddenly change direction.  Simulations are carried out in Matlab.\\

\begin{figure}[!htbp]
\begin{center} 
\includegraphics[width=5cm]{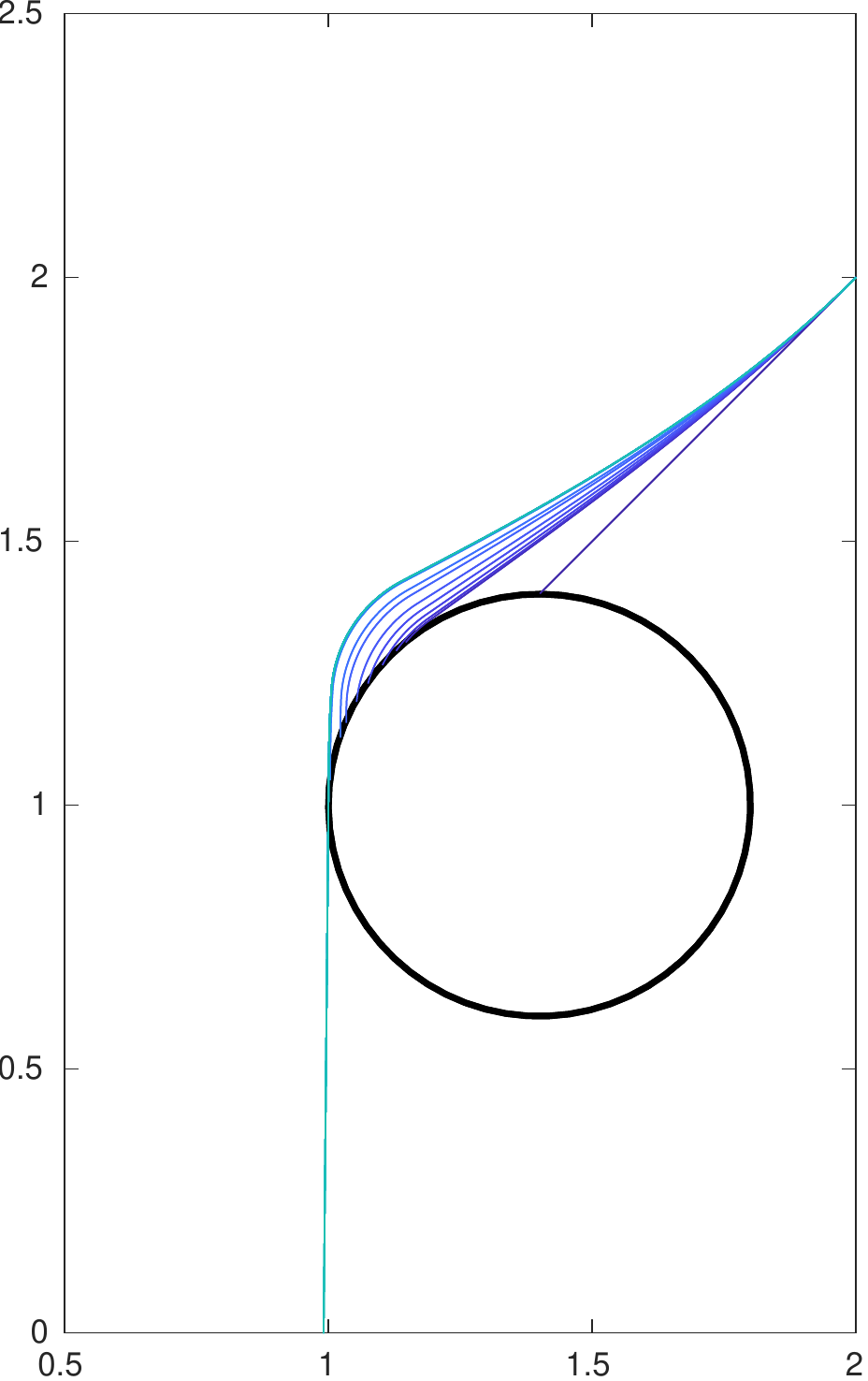}
\caption{\small  Artificial root not in breakdown configuration whose goal is to growth in depth}\label{sim1}
\end{center}
\end{figure}

\textbf{Simulation 1} In the first simulation the obstacle is a disc with center $(a,b)=(1.4,1)$ and radius $r=0.3$. The root has origin in $(2,2)$ and initial shape $y=x$ for $1.4\leq x \leq 2$. We fix the bound on the curvature $\kappa_0=4$. As we can see in Fig $\ref{sim1}$, the root bends, avoid the obstacle and grows vertically downward.  

\begin{figure}[!ht]
\begin{center}
\includegraphics[width=6cm]{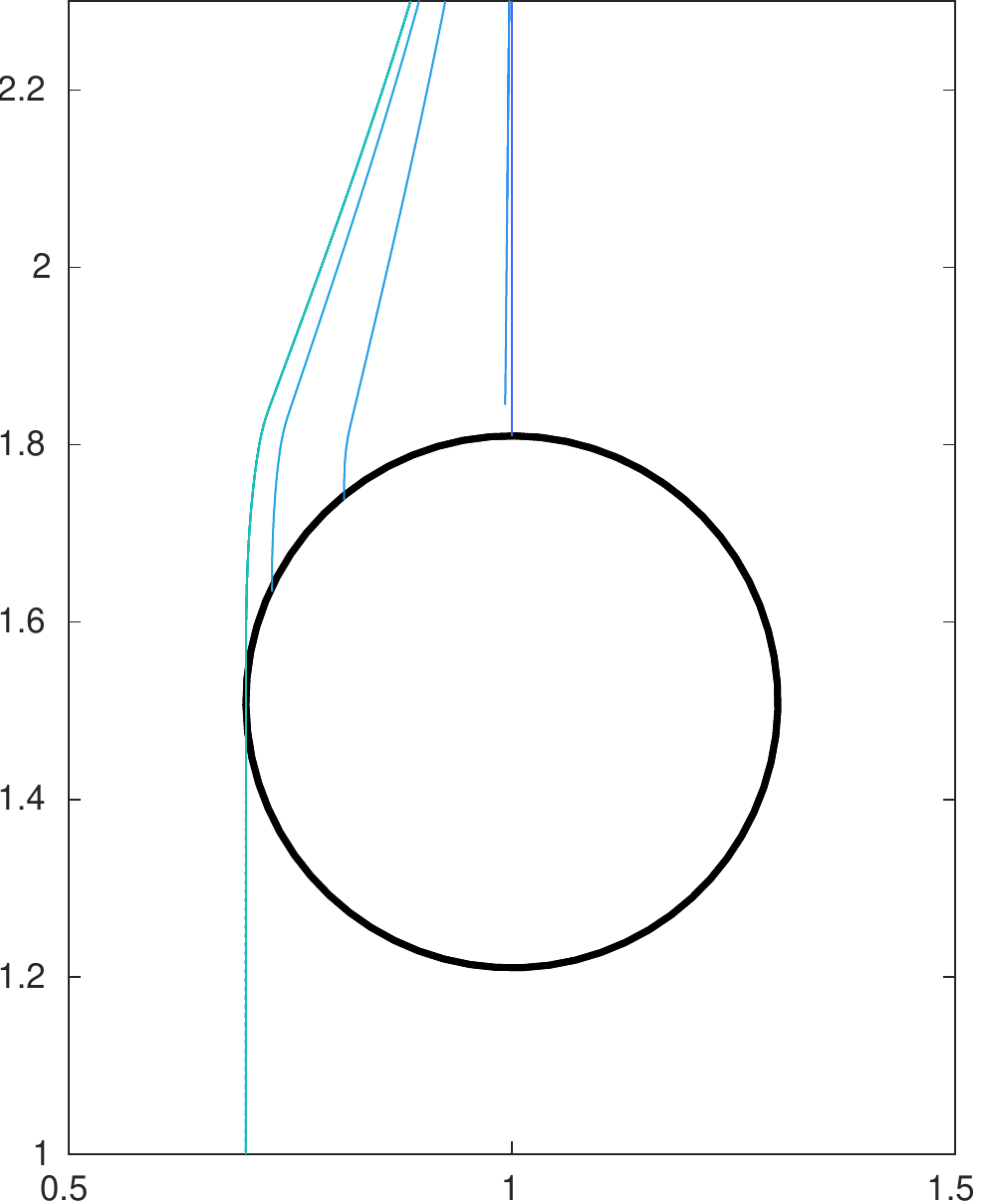}
\quad\quad
\includegraphics[width=6cm]{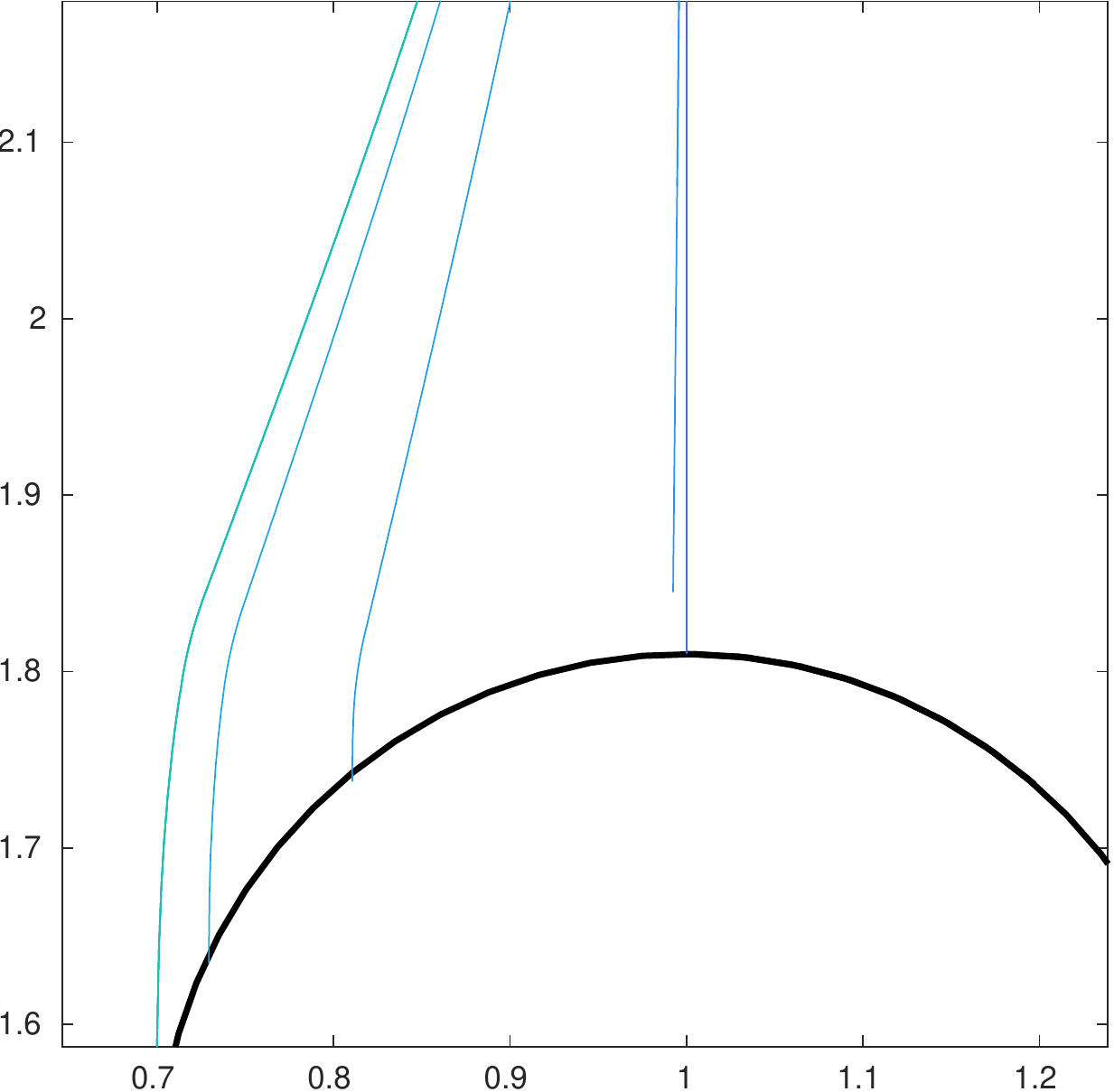}
\caption{\small  Artificial root in breakdown configuration which restarts }\label{sim2}
\end{center}
\end{figure}

\textbf{Simulation 2} Here we want to show how the restarting procedure works when the root is in a breakdown configuration. The obstacle is now a disc with center $(a',b')=(1,1.5)$ and radius $r=0.3$. The root has origin in $(1,2.5)$ and initial shape is $x=1$ for $1.8\leq y \leq 2.5$. The bound on the curvature is again $\kappa_0=4$. As we can see in Fig $\ref{sim2}$, the original root meets the obstacle perpendicularly with no curvature, hence by applying the restarting algorithm $(\ref{R1})$ the new attempt  grows far from the previous one.  In this way we can  avoid the region already explored.\\

Compared with the simulations performed in \cite{MR3650339},
the main difference here is that the evolution of our curve is now described by 
two different equations, one  for the body and one with the control 
 for the tip. Thus, the procedure for adding
nodes is different, since it is also determined by the feedback control. \\
Furthermore the robotic root must be able to recoil and restart whenever needed.
This requires to take into account all previous attempts, in order not to travel 
close to trajectories already investigated.
For simplicity, in the above simulations we have neglected the additional function $h$,  describing the hardness of the soil.
This corresponds to the case where the root descends into an empty cavity, between rigid obstacles.

\v

{\bf Acknowledgments.}  The authors would like to thank Prof. Piero Marcati for suggesting the problem.
The first author 
was partially supported by the Istituto Nazionale di Alta Matematica ``F.~Severi", through GNAMPA.
The research of the second and third authors was partially supported by 
NSF, with grant
DMS-1714237 ``Models of controlled biological growth".

\bibliographystyle{siam}
\bibliography{reference-root}

\end{document}